\documentclass[12pt,oneside]{article}

\usepackage{amssymb}

\usepackage{amsmath, amsthm}

\usepackage{amsfonts}

\usepackage{longtable}

\newtheorem{theorem}{Theorem}[section]
\newtheorem{lemma}[theorem]{Lemma}

\newtheorem{proposition}[theorem]{Proposition}

\newtheorem{corollary}[theorem]{Corollary}

\newtheorem{definition}[theorem]{Definition}

\newenvironment{ack}
{\begin{trivlist}  \item \textsc{Acknowledgments.}~}
{\end{trivlist}}

\newenvironment{proof of claim}
{\begin{trivlist}  \item \textsc{Proof of Claim:}~} {\hfill $\Box$
\end{trivlist}}

\newcommand{\defn}[1]{\ensuremath{\mathrm{Def^n}}(#1)}

\newcommand{\closure}[1]{\ensuremath{\mathrm{cl}}(#1)}
\newcommand{\interior}[1]{\ensuremath{\mathrm{int}}(#1)}

\def \st {\operatorname{st}}

\def \R{\mathbb{R}}
\def \Q{\mathbb{Q}}
\def \N{\mathbb{N}}
\def \m{\mathfrak{m}}
\def \k{\boldsymbol{k}}
\def \inc{\operatorname{i}}
\def \dec{\operatorname{d}}

\def \ind{\operatorname{ind}}

\newcommand{\ma}{\mathfrak{m}}

\def \Def{\operatorname{Def}}

\begin{document}

\author{Jana Ma\v{r}\'{i}kov\'{a}\\
Dept. of math. and stat., McMaster
University\\marikova@math.mcmaster.ca}
\title{O-minimal fields with standard part map}

\maketitle

\begin{abstract}
Let $R$ be an o-minimal field and $V$ a proper convex subring with
residue field $\boldsymbol{k}$ and standard part (residue) map $\st
\colon V\to \boldsymbol{k}$.  Let $\boldsymbol{k}_{\ind}$ be the
expansion of $\boldsymbol{k}$ by the standard parts of the definable
relations in $R$. We investigate the definable sets in
$\boldsymbol{k}_{\ind}$ and conditions on $(R,V)$ which imply
o-minimality of $\boldsymbol{k}_{\ind}$. We also show that if $R$ is
$\omega$-saturated and $V$ is the convex hull of $\Q$ in $R$, then
the sets definable in $\boldsymbol{k}_{\ind}$ are exactly the
standard parts of the sets definable in $(R,V)$.
\end{abstract}

\begin{section}{Introduction}

Throughout $R$ is an o-minimal field, that is, an o-minimal
expansion of a real closed field, and $V$ is a proper convex subring
with maximal ideal $\ma$, ordered residue field $\boldsymbol{k}=V/\ma$,
and standard part (residue) map $\st \colon V \to \boldsymbol{k}$.
This map induces a map $\st \colon V^n \to \boldsymbol{k}^n$ and
for $X\subseteq R^n$ we put $\st X:= \st (X\cap V^n )$.  By
$\boldsymbol{k}_{\ind}$ we denote the ordered field $\boldsymbol{k}$
expanded by the relations $\st X$ with $X\in \defn{R}$, $n=1,2,\dots$.
Unless indicated otherwise, by ``definable" we mean
``definable with parameters in the structure $R$".

The most important case of a convex subring of $R$ is the convex
hull
$$\mathcal{O}:= \{x\in R:\ |x|\le q\text{ for some }q\in
\mathbb{Q}^{>0}\}$$
of $\mathbb{Q}$ in $R$. If $V=\mathcal{O}$, then the ordered field
$\boldsymbol{k}$ is archimedean and we identify $\boldsymbol{k}$
with its image in the ordered field $\R$ of real numbers via the
unique ordered field embedding of $\boldsymbol{k}$ into $\R$. In
particular, if $R$ is $\omega$-saturated and $V=\mathcal{O}$, then
$\boldsymbol{k}=\R$.

\medskip\noindent
We consider the following questions:
\begin{enumerate}
\item[(1)] Under what conditions on $(R,V)$ is $\boldsymbol{k}_{\ind}$
o-minimal?
\item[(2)] How complicated are the definable relations of
$\boldsymbol{k}_{\ind}$ in terms of the basic relations $\st X$ with
definable $X\subseteq R^n$?
\end{enumerate}
Here is a brief history of these problems.  In 1983, Cherlin and
Dickmann \cite{chd} proved quantifier elimination for real closed
fields with a proper convex subring.  In 1995 van den Dries and
Lewenberg \cite{tconv} identified the notion of {\em $T$-convex
subring\/} of an o-minimal field as a suitable analogue of {\em
convex subring of a real closed field\/} (here $T$ is the theory of
the given o-minimal field). A convex subring $V$ of $R$ is said to
be $\text{Th}(R)$-convex if $f(V)\subseteq V$ for every continuous
$\emptyset$-definable function $f \colon R\to R$. The situation when
$V$ is a $\text{Th}(R)$-convex subring of $R$ is well-understood;
see \cite{tconv} and \cite{tconvex}. In particular,
$\boldsymbol{k}_{\ind}$ is o-minimal in that case.


\medskip\noindent The structure $\boldsymbol{k}_{\ind}$ is not always
o-minimal, as the example on page \pageref{example} shows.  A
theorem by Baisalov and Poizat \cite{bp} implies that
$\boldsymbol{k}_{\ind}$ is always weakly o-minimal. Hrushovski,
Peterzil and Pillay observe in \cite{nip} that if $R$ is
sufficiently saturated and $V=\mathcal{O}$, then it follows from
\cite{bp} that $\boldsymbol{k}_{\ind}$ is o-minimal, because then
$\boldsymbol{k}=\R$ and for expansions of the ordered field $\R$
weak o-minimality is the same as o-minimality.
However, \cite{nip} gives no information about question (2) in that
situation, which includes cases where $\mathcal{O}$ is not
$\text{Th}(R)$-convex; we say more about this in the remark on page
\pageref{pseudoreal}.

\medskip\noindent {\bf Good cell decomposition. \/}
In \cite{st} we answered (2) for the situation in \cite{nip} by
means of {\em good cell decomposition}, which also gives the
o-minimality of $\R_{\ind}$ without using \cite{bp}. In
the present paper we obtain good cell decomposition (and thus
o-minimality of $\boldsymbol{k}_{\ind}$) under more general {\em
first-order\/} assumptions on the pair $(R,V)$. More precisely,
suppose $(R,V)\models\Sigma_{\inc}$ where $\Sigma_{\inc}$ is defined below. Theorem
2.21 says that then the subsets of $\boldsymbol{k}^n$ definable in
$\boldsymbol{k}_{\ind}$ are the finite unions of differences $\st X
\setminus \st Y$, where $X, Y \subseteq R^n$ are definable. It
follows that $\boldsymbol{k}_{\ind}$ is o-minimal. Theorem 2.21 is
proved in the same way as the corresponding theorem in \cite{st},
except that uses of saturation in \cite{st} are replaced by uses of
$\Sigma_{\inc}$. Also the proof of Lemma 4.1. in \cite{st} does not
generalize to our setting, and this is replaced here by a more
elementary proof of Lemma \ref{nice function} below.


\bigskip\noindent
The following conditions on $(R,V)$ are related to good cell
decomposition. To state these, let $I:= \{x\in R:\ |x| \le 1\}$, and
for $X\subseteq R^{1+n}$ and $r\in R$, put $$X(r):=\{ x\in R^n :\;
(r,x)\in X \}.$$ We let $\ma^{>r} := \{x \in \ma : \; x>r\}$ for
$r\in \ma$. We define the conditions $\mathcal{I}$, $\Sigma_{\inc}$,
$\Sigma_{\dec}$, $\Sigma$, and $\mathcal{C}$ on pairs $(R,V)$ as
follows:
\begin{enumerate}\label{conditions}
\item[($\mathcal{I}$)] if $X,Y \subseteq I^n$ are definable, then there
is a definable $Z\subseteq I^n$ such that $\st X\ \cap\ \st Y\ =\
\st Z$;
\item[($\Sigma_{\inc}$)] if $X\subseteq I^{1+n}$ is definable
and $X(r)\subseteq X(s)$ for all $r,s\in I$ with $r \le s$, then
there is $\epsilon_0 \in \ma^{>0}$ such that $\st X(\epsilon_0 )
=\st X(\epsilon )$ for all $\epsilon\in \ma^{>\epsilon_0 }$;
\item[($\Sigma_{\dec}$)] if $X\subseteq I^{1+n}$ is definable
and $X(r)\supseteq X(s)$ for all $r,s\in I$ with $r \le s$, then
there is $\epsilon_0 \in \ma^{>0}$ such that $\st X(\epsilon_0 )
=\st X(\epsilon )$ for all $\epsilon\in \ma^{>\epsilon_0 }$;
\item[($\Sigma$)] if $X\subseteq I^{1+n}$ is definable, then
there is $\epsilon_0 \in \ma^{>0}$ such that $\st X(\epsilon_0 )
=\st X(\epsilon )$ for all $\epsilon\in \ma^{>\epsilon_0 }$;
\item[($\mathcal{C}$)] the $\boldsymbol{k}_{\ind}$-definable closed
subsets of $\boldsymbol{k}^{n}$ are exactly the sets $\st X$ with
definable $X\subseteq R^n$.
\end{enumerate}
One should add here ``for all $n$ and $X,Y$'' as initial clause to
$\mathcal{I}$, and likewise with the other conditions.  In Section 3
we prove that for all $(R,V)$,
\begin{trivlist}
\item[a)] \qquad \qquad $\mathcal{I}\ \Longleftrightarrow \Sigma_{\inc}$;
\item[b)] \qquad \qquad $\Sigma_{\inc} \ \Longrightarrow \ \boldsymbol{k}_{\ind}$ is o-minimal;
\item[c)] \qquad \qquad $\Sigma\ \Longrightarrow\ \mathcal{C}$.
\end{trivlist}
We do not know whether the converse of b) holds.  In a subsequent
paper with van den Dries \cite{jalou} we shall prove the converse of
c), and also $\Sigma_{\inc} \Longrightarrow \mathcal{C}$, yielding
$\Sigma_{\inc} \Longleftrightarrow \Sigma$.

Our definition of $\mathcal{I}$ is not of first-order nature, but by
a) it is equivalent to first-order conditions.  Similarly
$\mathcal{C}$ will turn out to be equivalent to first order
conditions by c) and its converse in \cite{jalou}.

\medskip\noindent
In Section 3 we also show that $(R,V)$ satisfies $\Sigma$ if any of
the following holds:
\begin{enumerate}

\item[(i)] $\operatorname{cofinality}(\ma)>2^{|\boldsymbol{k}|}$;
\item[(ii)] $V$ is $T$-convex, where $T:= \text{Th}(R)$;
\item[(iii)] $R$ is $\omega$-saturated and $V=\mathcal{O}$.
\end{enumerate}

\medskip\noindent {\bf Traces. \/} Call a
set $X\subseteq R^n$ a {\em trace\/} if $X=Y\cap R^n$ for some
definable $n$-ary relation $Y$ in some elementary extension of $R$,
where we allow parameters from that elementary extension to define
$Y$.  In Section 4 we assume that $R$ is $\omega$-saturated and
$V=\mathcal{O}$, and under these assumptions we
characterize the definable sets in $\R_{\ind}$ in terms of
traces. As a corollary we obtain that if $R$ is $\omega$-saturated
and $V=\mathcal{O}$, then
$$\Def^n(\R_{\ind})\ =\ \{\st X:\ X\in \Def^n(R,\mathcal{O})\}.$$

\medskip\noindent
We do not know if the analogue of this corollary holds under the
more general first-order assumption $\Sigma$. We do know that if $V$
is $\text{Th}(R)$-convex, then, for all $n$,
$$\Def^n (\boldsymbol{k}_{\ind }) = \{ \st X: X\in \Def^n (R,V)
\}.$$

\medskip\noindent {\bf Remark.\/}
\label{pseudoreal} In 1996 van den Dries \cite{lou-qu} asked the
following question: Let $L$ be a language extending the language of
ordered rings, and let $T(L,\R )$ be the set of all sentences true
in all $L$-expansions of the real field.  Call $R$ {\em
pseudo-real\/} if $R\models T(L,\R )$. Is every o-minimal field
pseudo-real?

If $R$ has an archimedean model, then $R$ is pseudo-real,
but the converse fails. Consider for example a proper
elementary extension of the real field and extend its language by a
name for an element $\lambda > \R$.  Then the theory of $R$ in the
extended language does not have an archimedean model but $R$ is of
course pseudo-real as a structure for this extended language.

In 2006 Lipshitz and Robinson \cite{lr} considered the ordered
Hahn field $\R ((t^{\Q}))$ with operations given by overconvergent
power series, and they proved its o-minimality. In 2007 Hrushovski
and Peterzil \cite{uk} showed that this Lipshitz-Robinson field is
not pseudo-real.
It is easy to see that if $R$ is a model of the
theory $T$ of the Lipshitz-Robinson field, then $\mathcal{O} \subseteq
R$ is not $T$-convex.

\medskip\noindent {\bf Preliminaries.\/}
We assume familiarity with o-minimal structures and their basic
properties; see for example \cite{lou}. Throughout we let $m,n$
range over the set $\N=\{0,1,2,\dots\}$ of natural numbers. Given a
one-sorted structure $\mathcal{M}=(M;\cdots)$ we let
$\Def^n(\mathcal{M})$ be the boolean algebra of definable subsets of
$M^n$. Let $K$ be an ordered field. For $x\in K$ we put $|x|:=
\max\{x,-x\}$, for $a=(a_1,\dots, a_n)\in K^n$ we put
$$|a|:= \max \{|a_i|:\ i=1,\dots,n\}\ \text{ if }n>0, \quad
|a|:=0 \text{ if }n=0,$$ and for $a,b\in K^n$ we put $d(a,b):=
|a-b|$. A {\em box in $K^n$\/} is a cartesian product of open
intervals
  $$(a_1 -\delta , a_1 +\delta )\times \dots \times (a_n -\delta
, a_n + \delta ),$$ where $a=(a_1,\dots,a_n) \in K^n$ and $\delta
\in K^{>0}$. A {\em $V$-box in $R^n$\/}\label{defVbox} is a box in
$R^n$ as above where $a\in V^n$ and $\delta \in V^{>\ma}$. So if
$B\subseteq R^n$ is a $V$-box, then $B\subseteq V^n$ and $\st B$
contains a box in $\boldsymbol{k}^n$.

An {\em interval\/} is always a nonempty open interval $(a,b)$ in
$R$, or in $\R$, or in $\boldsymbol{k}$, as specified. We already
defined $I:= \{x\in R:\ |x|\le 1\}$ and more generally, for each
ordered field $K$ we put $I(K):= \{x\in K:\ |x|\le 1\}$. For $a\in
R^n$ and definable nonempty $X\subseteq R^n$ we set
$$d(a,X):=\inf\{d(a,x):\; x\in X \},$$
and likewise for $a \in \boldsymbol{k}^n$ and definable nonempty $X
\subseteq \boldsymbol{k}^n$ when $\boldsymbol{k}_{\ind}$ is
o-minimal. A set $X\subseteq R^n$ is said to be {\em  $V$-bounded\/}
if there is $a\in V^{>0}$ such that $|x|\le a$ for all $x\in X$.
(For $V=\mathcal{O}$ this is the same as {\em strongly bounded}.)
The {\em hull\/} of $X\subseteq \boldsymbol{k}^{n}$ is the set $X^h
:= \st^{-1}(X)\subseteq V^n$.

Given sets $X,Y$ and $S\subseteq X\times Y$ we put
$$S(x):= \{y\in Y:\ (x,y)\in S\}.$$
If $X$ is a subset of an ambient set $M$ that is understood from the
context, then
$$X^c := \{x\in M :\; x\not\in X \}.$$
We often use the following projection maps for $m\leq n$:
\begin{align*}
p^{n}_{m} &: R^n \rightarrow R^m , \qquad(x_1 , \dots ,x_n ) \mapsto
(x_1 , \dots ,x_m )\\
\pi^{n}_{m}&: \boldsymbol{k}^n \rightarrow \boldsymbol{k}^m, \qquad (x_1
, \dots , x_n ) \mapsto (x_1 , \dots ,x_m ).
\end{align*}
 Given a map $f \colon X \to Y$ we let
$$ \Gamma f := \{(x,y)\in X\times Y:\ f(x)=y\}$$
denote its graph.
\end{section}
\begin{ack}
This paper contains some of the results in the author's PhD thesis.
The author would like to thank her advisor Lou van den Dries for
advice.
\end{ack}

\begin{section}{Good cell decomposition}

\begin{subsection}{General facts on standard part sets} Recall that $R$ is an
o-minimal field and $V$ is a proper convex subring of $R$. We begin
with some results requiring no extra assumption on $(R,V)$. A very
useful fact of this kind is the $V$-box Lemma (Corollary
\ref{vboxlemma}).

\begin{lemma}\label{st-are-closed}
If $X \subseteq R^n$ is definable, then $\st X$ is closed.
\end{lemma}
\begin{proof}
Let $X \subseteq R^n$ be definable and assume towards a
contradiction that we have an $a \in \closure{\st X } \setminus \st
X$. Take $a' \in R^n$ such that $\st{a'}=a$. Then, by o-minimality
of $R$, $d(a' ,X)$ exists in $R$ and $d(a' ,X) >
\frak{m}$. So there is a neighborhood $U \subseteq \boldsymbol{k}^n$
of $a$ with $U \cap \st X = \emptyset$, a contradiction.
\end{proof}

Let $\text{St}_n$ be the collection of all sets $\st X$ with
definable $X\subseteq R^n$. Note that if $X,Y\in \text{St}_n$, then
$X\cup Y\in \text{St}_n$; if $X\in \text{St}_m$ and $Y\in
\text{St}_n$, then $X\times Y\in \text{St}_{m+n}$. The next lemma is
almost obvious. To state it we use the projection maps
$\pi=\pi^{m+n}_m \colon \boldsymbol{k}^{m+n}\to \boldsymbol{k}^m$
and $p=p^{m+n}_m \colon R^{m+n} \to R^m$.

\begin{lemma}\label{stan0} Let $X\in \operatorname{St}_{m+n}$. Then
\begin{enumerate}
\item[$(1)$] if $X$ is bounded, then $\pi(X)\in
\operatorname{St}_m$;
\item[$(2)$] if $X=\st X'$ where the set $X'\subseteq R^{m+n}$ is
definable in $R$ and satisfies $X'\cap p^{-1}(V^m)\subseteq
V^{m+n}$, then $\pi(X)\in \operatorname{St}_m$.
\end{enumerate}
\end{lemma}

\begin{lemma}\label{stan2} If $X\subseteq R$ is definable, then
$\st X$ is a finite union of intervals and points in
$\boldsymbol{k}$.
\end{lemma}
\begin{proof} This is immediate from the o-minimality of $R$.
\end{proof}

\noindent Recall the definition of a $V$-box from page
\pageref{defVbox}.
Below $p$ is the projection map $R^{n+1} \rightarrow R^{n}$ given by
$p(x_1 , \dots ,x_{n+1})= (x_1, \dots ,x_n )$.

\begin{lemma}\label{nice function}
\indent
\begin{enumerate}
\item[{\em ($A_n$)}] If $D \subseteq V^{n+1}$ is a $V$-box,
and $f \colon Y \rightarrow R$, where $Y \subseteq V^{n}$, is
definable and continuous with $f(Y)\subseteq V$, then there is a
$V$-box $B \subseteq D$ with $B \cap \Gamma f=\emptyset$.
\item[{\em ($B_n$)}] If $D \subseteq V^n$ is a $V$-box, and
$\mathcal{C}$ is a decomposition of $D$, then there is $C \in
\mathcal{C}$ such that $C$ contains a $V$-box.
\end{enumerate}
\end{lemma}
\begin{proof}
It is clear that $(B_1 )$ holds. We first show that $(B_n )$ implies
$(A_n )$. Let $f \colon Y \rightarrow V$ be definable and
continuous, with $Y \subseteq V^n$, and let
$$D = (a_1 , b_1 ) \times \dots \times (a_{n+1} , b_{n+1} )
\subseteq V^{n+1}$$ be a $V$-box. Take $p , q \in V$ such that
$a_{n+1} <p < q <b_{n+1}$ and $$q-p , p-a_{n+1} , b_{n+1} -q >
\frak{m},$$ and pick $\delta
> \frak{m}$ with $\delta < \min \{p - a_{n+1} , \frac{q - p}{2} ,
b_{n+1} -q \}$. Define
\begin{align*}
X(p) &:= \{x \in p^{n+1}_{n} D \cap Y:\;
f(x) \in (p-\delta , p + \delta )    \}\\
X(q) &:= \{x \in p^{n+1}_{n} D \cap Y:\;  f(x) \in (q - \delta , q +
\delta ) \},
\end{align*}
and note that $X(p) \cap X(q)=\emptyset$. Take a decomposition
$\mathcal{C}$ of $R^{n}$ such that $\mathcal{C}$ partitions the sets
$p^{n+1}_{n} D$, $X(p)$, and $X(q)$. By $(B_n )$, there is $C \in
\mathcal{C}$ such that $C \subseteq p^{n+1}_{n} D$ and $C$ contains
a $V$-box $P$. Then $P \times (p-\delta , p +\delta )$ or $P \times
(q-\delta , q+\delta )$ yields the desired $V$-box $B$.

Next, we show that $(A_n )$ and $(B_n )$ imply $(B_{n+1} )$. Let $D
\subseteq V^{n+1}$ be a $V$-box and let $\mathcal{C}$ be a
decomposition of $D$. Then $p^{n+1}_{n} \mathcal{C}$ is a
decomposition of $p^{n+1}_{n} D$ and by $(B_n )$ we can take $C \in
\mathcal{C}$ such that $p^{n+1}_{n} C$ contains a $V$-box $P$. Let
$C_1, \dots ,C_k$ be the cells in $\mathcal{C}$ such that
$p^{n+1}_{n} C = p^{n+1}_{n} C_i$ for $i=1, \dots ,k$. After
restricting the functions $p^{n+1}_{n} C \rightarrow R$ used to
define $C_1 , \dots , C_k$ to $P$ we see that it is enough to prove
the following:

Let $f_1 , \dots ,f_m \colon P \rightarrow V$ be definable and
continuous and let
$p,q \in V$ be such that $p<q$ and $|q-p|
>\ma$. Then there is a $V$-box $B \subseteq P \times (p,q)$
with $B \cap \Gamma f_j = \emptyset$ for all $j$.

For $m=1$ this statement follows from ($A_n$), and for $m>1$ it
follows by a straightforward induction on $m$ using again ($A_n$).
\end{proof}

\begin{corollary}{\bf ($V$-Box Lemma)}\label{vboxlemma}
Let $X \subseteq R^n$ be definable and let $D\subseteq
\boldsymbol{k}^n$ be a box such that $D \subseteq \st X$. Then $X$
contains a $V$-box $B$ with $\st B \subseteq D$.
\end{corollary}

\begin{proof} We may assume that $X \subseteq V^n$, and that
$\closure{D} \subseteq \st X$. Pick a $V$-box $D'\subseteq R^n$ such
that $\st D'=\closure{D}$, and take a decomposition $\mathcal{C}$ of
$R^n$ which partitions both $D'$ and $X$.
By Lemma \ref{nice function}, we can take $C \in \mathcal{C}$ such
that $C \subseteq D'$ and $C$ contains a $V$-box $B$. It is clear
that $B\cap X \not= \emptyset$, otherwise $D$ would contain a box
whose intersection with $\st X$ is empty. So $B\subseteq C\subseteq
X$.
\end{proof}

\begin{corollary}
If $X \subseteq R^n$ is definable, then $\st (X) \cap \st (X^{c})$
has empty interior in $\boldsymbol{k}^n$.
\end{corollary}
By \cite{bp}, $\boldsymbol{k}_{\ind}$ is weakly o-minimal.
MacPherson, Marker and Steinhorn define in \cite{womin} a notion of
dimension for weakly o-minimal structures:
\begin{definition}
Let $M$ be a weakly o-minimal structure, and let $X\subseteq M^n$ be
definable in $M$. If $X\not=\emptyset$, then $\dim_{w}(X)$ is the
largest integer $k \in \{ 0,\dots ,n  \}$ for which there is a
projection map
$$p \colon M^n \rightarrow M^k, \qquad(x_1 , \dots ,x_n )\mapsto
(x_{\lambda(1)}, \dots , x_{\lambda (k)}),$$ where $1\leq \lambda
(1)<\dots <\lambda (k)\leq n$, such that $\interior{pX}\not=
\emptyset$. We set $\dim_{w} (\emptyset )=-\infty$.
\end{definition} Note that if $M$ is o-minimal, then the above
notion of dimension agrees with the usual dimension for o-minimal
structures.
\begin{corollary}
$\dim_{w}{(\st X)} \leq \dim{(X)}$ for $V$-bounded $X\in \Def^n
(R)$.
\end{corollary}

\end{subsection}

\begin{subsection}{Good cells}
We define good cells in analogy with \cite{st}, and we state some
results needed in the proof of good cell decomposition. We omit
proofs that are as in \cite{st}.

\begin{definition}
Given functions $f\colon X \rightarrow R$ with $X \subseteq R^n$,
and $g\colon C \rightarrow \boldsymbol{k}$ with $C \subseteq
\boldsymbol{k}^n$, we say that $f$ {\em induces\/} $g$ if $f$ is
definable $($so $X$ is definable$)$, $C^h\subseteq X$, $f|C^h$ is
continuous, $f(C^{h}) \subseteq V$ and $\Gamma g = \st(\Gamma f)
\cap (C\times \boldsymbol{k})$.
\end{definition}

\begin{lemma}
Let $C \subseteq \boldsymbol{k}^n$ and suppose $g\colon  C
\rightarrow \boldsymbol{k}$ is induced by the function $f\colon X
\rightarrow R$ with $X \subseteq R^n$. Then $g$ is continuous.
\end{lemma}
\begin{proof}
Assume towards a contradiction that $g$ is not continuous at $c \in
C$. Let $r \in \boldsymbol{k}^{>0}$ be such that for every
neighborhood $B \subseteq \boldsymbol{k}^n$ of $c$ there is $b \in
B\cap C$ with $|g(c) - g(b)| \geq r$. Pick $c' \in R^n$ with
$\st{c'}=c$ and define
$$Y:=\{ x \in X:\; |f(c') - f(x)| \geq \frac{r'}{2} \},$$ where $r'
\in R^{>0}$ is such that $\st{r'}=r$. Then $d(c',Y)$ exists
in $R$. If $d(c',Y)$ is infinitesimal then, since $Y$ is
closed, there is $y \in Y$ such that $\st{y} = \st{c'}$, a
contradiction with $f$ inducing a function. Hence $d(c',Y)
> \frak{m}$, but this yields a neighborhood $B \subseteq
\boldsymbol{k}^n$ of $c$ such that $g(B \cap C) \subseteq (g(c)-r ,
g(c)+r)$, a contradiction.
\end{proof}

\noindent For $C\subseteq \boldsymbol{k}^n$ we let $G(C)$ be the set
of all $g\colon C \to \boldsymbol{k}$ that are induced by some
definable $f\colon X \to R$ with $X\subseteq R^n$.

\begin{lemma}\label{stan4}
Let $1\le j(1) < \dots < j(m)\le n$ and define $\pi \colon
\boldsymbol{k}^n \to \boldsymbol{k}^m$ by
$$\pi(x_1,\dots,x_n)=
(x_{j(1)}, \dots,x_{j(m)}).$$ Let $C\subseteq \boldsymbol{k}^n$ and
suppose $g\in G(\pi C)$. Then $g\circ \pi|_{C} \in G(C)$.
\end{lemma}

\begin{definition}\label{defgood} Let $i=(i_1,\dots,i_n)$ be a sequence of
zeros and ones. Good $i$-cells are subsets of $\boldsymbol{k}^n$
obtained by recursion on $n$ as follows:
\begin{enumerate}
\item[{\em (i)}] For $n=0$ and $i$ the empty sequence,
the set $\boldsymbol{k}^0$ is the only good $i$-cell, and for $n=1$,
a {\em good $(0)$-cell\/} is a singleton $\{a\}$ with $a\in
\boldsymbol{k}$; a {\em good $(1)$-cell\/} is an interval in
$\boldsymbol{k}$.

\item[{\em (ii)}] Let $n>0$ and assume inductively that good $i$-cells are
subsets of $\boldsymbol{k}^n$. A good $(i,0)$-cell is a set $\Gamma
h \subseteq \boldsymbol{k}^{n+1}$ where $h \in G(C)$ and $C
\subseteq \boldsymbol{k}^n$ is a good $i$-cell. A good $(i,1)$-cell
is either a set $C\times \boldsymbol{k}$, or a set $(-\infty,
f)\subseteq \boldsymbol{k}^{n+1}$, or a set $(g,h) \subseteq
\boldsymbol{k}^{n+1}$, or a set $(f, +\infty) \subseteq
\boldsymbol{k}^{n+1}$, where $f,g,h \in G(C)$, $g<h$, and $C$ is a
good $i$-cell.
\end{enumerate}
\end{definition}

\noindent One verifies easily that a good $i$-cell is open in
$\boldsymbol{k}^n$ iff $i_1=\dots =i_n=1$, and that if
$i_1=\dots=i_n=1$, then every good $i$-cell is homeomorphic to
$\boldsymbol{k}^n$. A {\em good cell in $\boldsymbol{k}^n$} is a
good $i$-cell for some sequence $i=(i_1,\dots,i_n)$ of zeros and
ones.

\begin{lemma}\label{inverse}
Let $C\subseteq \boldsymbol{k}^n$ be a good $(i_1,\dots,i_n)$-cell,
and let $k\in \{1,\dots,n\}$ be such that $i_k =0$. Let $\pi \colon
\boldsymbol{k}^n \to \boldsymbol{k}^{n-1}$ be given by
$$ \pi(x_1,\dots,x_n)=(x_1,\dots,x_{k-1}, x_{k+1}, \dots,x_n).$$
Then $\pi(C)\subseteq \boldsymbol{k}^{n-1}$ is a good cell, $\pi|C
\colon C \to \pi(C)$ is a homeomorphism, and if $E\subseteq \pi(C)$
is a good cell, so is its inverse image $\pi^{-1}(E)\cap C$.
\end{lemma}

\end{subsection}

\begin{subsection}{More on good cells} Recall the conditions $\mathcal{I}$ and
$\Sigma_{\inc}$ on pairs $(R,V)$ from page \pageref{conditions}. We
prove here that $(R,V)\models \mathcal{I}$
 iff $(R,V)\models \Sigma_{\inc}$. This yields that if $(R,V)\models \Sigma_{\inc}$, then good cells in $\boldsymbol{k}^n$ are
differences of standard parts of definable subsets of $R^n$.


It is not difficult to show that if $(R,V)\models \mathcal{I}$, then
for all $n$ and all definable $X,Y \subseteq R^n$ there is a
definable $Z\subseteq R^n$ such that $\st (X) \cap \st (Y) = \st Z$:
Set $J(\boldsymbol{k} ):= (-1,1) \subseteq \boldsymbol{k}$ and $J:=
(-1,1) \subseteq R$. We shall use the definable homeomorphism

$$\tau_n \colon R^n \rightarrow J^{n}\colon (x_1 \dots ,x_n ) \mapsto
(\frac{x_1 }{\sqrt{1 + x_{1}^2}} , \dots , \frac{x_n
}{\sqrt{1+x_{n}^{2}}}),$$ and we also let $\tau_n$ denote the
homeomorphism

$$\tau_n \colon \boldsymbol{k}^n \rightarrow J(\boldsymbol{k}
)^{n}\colon (x_1 \dots ,x_n ) \mapsto (\frac{x_1 }{\sqrt{1 +
x_{1}^2}} , \dots , \frac{x_n }{\sqrt{1+x_{n}^{2}}}).$$ One easily
checks that $\tau_1 \colon R \rightarrow J$ induces $\tau_1 \colon
\boldsymbol{k} \rightarrow J(\boldsymbol{k} )$, and that for $X\in
\Def^n (R)$,
$$\tau_n (\st X) = \st (\tau_n X )\cap J(\boldsymbol{k})^n \;\mbox{
and }\; \tau_{n}^{-1}(\st (X) \cap J(\boldsymbol{k})^n )= \st
(\tau_{n}^{-1}(X)),$$ where $\tau_{n}^{-1}\colon J^n \rightarrow
R^n$ and $\tau_{n}^{-1} \colon  J(\boldsymbol{k})^n \rightarrow
\boldsymbol{k}^n$ are the inverse functions of $\tau_n \colon R^n
\rightarrow J^n$ and of $\tau_n \colon \boldsymbol{k}^n \rightarrow
J(\boldsymbol{k})^n$ respectively.

\medskip\noindent
Suppose $(R,V)$ satisfies $\mathcal{I}$. Then for all $n$ and all
$X,Y\in \Def^n (R)$ there is $Z\in \Def^n (R)$ such that $\st
(X)\cap \st (Y)=\st (Z)$. To see this, let $X,Y \in \Def^n (R)$.
Then $\tau_n (X), \tau_n (Y)\subseteq J^n$, so we can take $Z\in
\Def^n (R)$ such that $$\st (\tau_n (X))\cap \st (\tau_n (Y))=\st
Z.$$ We claim that
$$\st (X)\cap \st (Y)= \st (\tau^{-1}_{n}(Z\cap J^n )).$$ To prove this it is enough to show that
\begin{equation}
\tau (\st (X)\cap \st (Y))= \tau(\st (\tau^{-1}_{n}(Z\cap J^n ))).
\end{equation}
Now the right-hand side of $(1)$ is equal to $$\st (Z\cap J^n ) \cap
J(\boldsymbol{k})^n = \st (Z)\cap J(\boldsymbol{k})^n ,$$ and we
have
$$\tau_n
(\st (X) \cap \st (Y))= \st (\tau_n X) \cap \st (\tau_n Y)\cap
J(\boldsymbol{k})^n .$$ In view of $\st (\tau_n (X))\cap \st (\tau_n
(Y))=\st Z$ this gives $(1)$.

In a similar way the condition $\Sigma_{\inc}$ implies its
``unrestricted version", i.e. the variant obtained by substituting
$R$ for $I$. We shall often use these facts silently.

\begin{lemma}\label{i implies sigma inc}
Suppose $(R,V)$ satisfies $\mathcal{I}$. Then $(R,V)\models
\Sigma_{\inc}$.
\end{lemma}
\begin{proof}
Let $X \subseteq I^{1+n}$ be definable and increasing in the first
variable. Towards proving that $X$ satisfies the conclusion of
$\Sigma_{\inc}$ we may assume that $X$ is closed.

\smallskip\noindent{\em Claim 1.} There is $\epsilon_0 \in \ma^{\geq
0}$ such that
$$\st (X) \cap (\{ 0 \} \times I(\boldsymbol{k})^n ) =
\st (X \cap ([0, \epsilon_0 ] \times I^n )).$$

\smallskip\noindent We set $Y:= \{ 0\} \times I^n$ and take a definable
$Z \subseteq I^{n+1}$ with $\st (X) \cap \st (Y) = \st (Z)$. We may
assume that $Z$ is closed and nonempty, and we set $\epsilon_1 :=
\sup\{d(z,X):\; z\in Z\}$ and $\epsilon_2 :=
\sup\{d(z,Y):\; z \in Z\}$. Then $\epsilon_1 , \epsilon_2 \in
\ma^{\geq 0}$, and we claim that $\epsilon_0 := \epsilon_1 +
\epsilon_2$ works. Clearly, $$\st (X \cap ([0,\epsilon_0 ]\times I^n
)) \subseteq \st (X) \cap (\{ 0 \}\times I(\boldsymbol{k})^n ).$$ So
let $a\in \st (X) \cap \st (Y)$. Then $a=\st z$ with $z\in Z$. We
have $d(z,X)\leq \epsilon_1 $ and $d(z,Y)\leq
\epsilon_2$. Since $Z$ is closed and $V$-bounded, we can take $x\in
X$ and $y\in Y$ such that $d(x,z)\leq \epsilon_1$,
$d(y,z)\leq \epsilon_2$. Then $d(x,y)\leq \epsilon_1 +
\epsilon_2 = \epsilon_0 $, and it follows that $$a = \st x \in \st
(X \cap ([0,\epsilon_0 ] \times I^n )).$$ This proves Claim 1. Let
$\epsilon_0$ be as in Claim 1.

\smallskip\noindent{\em Claim 2.} $\st X(\epsilon )=\st X(\epsilon_0
)$ for all $\epsilon \in \ma^{\geq \epsilon_0 }$.

\smallskip\noindent
It is clear that $\st X(\epsilon_0 )\subseteq \st X(\epsilon )$ for
all $\epsilon \geq \epsilon_0$. To prove the other inclusion, let $a
\in \st X(\epsilon )$ and take $x \in X(\epsilon )$ such that $\st x
= a$. Then $$(0,a) \in \st (X) \cap (\{ 0 \} \times
I(\boldsymbol{k})^n ),$$ hence $$(0,a) \in \st (X \cap
([0,\epsilon_0 ] \times I^n ))$$ by Claim 1. Because $X$ is
increasing in the first variable, this implies $(0,a) \in \st
X(\epsilon_0 )$.
\end{proof}

\begin{lemma}\label{stan1}
$\Sigma_{\inc}\ \Longrightarrow\ \mathcal{I}$.

\end{lemma}
\begin{proof} Suppose $(R,V)$ satisfies $\Sigma_{\inc}$. Let
$X,Y \subseteq I^n$ be definable and nonempty.  For $\epsilon \in
R^{\geq 0}$ define
$$Y^{\epsilon } := \{ x \in R^n : \;  d(x, Y) \leq \epsilon
\}.$$ We claim that $$\bigcup _{\epsilon }\st{ (X \cap Y^{\epsilon
})}= \st{X} \cap \st{Y},$$ where $\epsilon$ ranges over all positive
infinitesimals. If $a \in \st{(X \cap Y^{\epsilon })}$, then clearly
$a \in \st{X}$ and $a \in \st{Y}$. If $a \in \st{X} \cap \st{Y}$,
then we can take $a' \in X$ and $a'' \in Y$ such that $\st a' = \st
a'' = a$ and $d(a',a'')<\epsilon$ for some $\epsilon \in
\ma^{>0}$. Hence $a' \in X\cap Y^{\epsilon }$.

Now by $\Sigma_{\inc}$, there is a positive infinitesimal
$\epsilon_0$ such that
$$\st{(X \cap Y^{\epsilon_0 })} = \bigcup_{\epsilon }\st{(X \cap
Y^{\epsilon })}.$$
\end{proof}
The proofs of the following two lemmas are similar to the proofs of
their counterparts in \cite{st}.

\begin{lemma}\label{stan1a} Suppose $(R,V)$ satisfies $\mathcal{I}$, and
let $X\subseteq R^n$ and $f\colon  X \to R$ be definable, and put
$$  X^{-}:=\{x\in X:\ f(x) < V\}, \quad
X^{+}:=\{x\in X:\ f(x)> V\}.$$ Then $\st (X^{-})$ and $\st(X^{+})$
belong to $\operatorname{St}_n$.
\end{lemma}

\begin{corollary} If $(R,V)$ satisfies $\mathcal{I}$, and $X\subseteq R^n$
and $g\colon  X \to R$ are definable, then $\st(\{x\in X: g(x)\in
\ma \}) \in \operatorname{St}_n$.
\end{corollary}

Conversely, if the conclusion of this corollary holds for all $n$
and definable $g\colon  X \to R$ with $X\subseteq R^n$, then $(R,V)$
satisfies $\mathcal{I}$. To see this, let $X,Y\subseteq V^n$ be
definable with $Y\ne \emptyset$. Assume the conclusion of the
corollary holds for the function $x \mapsto d(x,Y)\colon  X
\to R$. Then we have a definable $Z\subseteq V^n$ such that $\st(Z)
= \st (\{ x \in X:\  d(x,Y)\in \ma \})$. This gives $\st (X)
\cap \st (Y)= \st (Z)$.

\medskip\noindent {\bf From now on until the end of Section 2 we
assume $(R,V)\models \Sigma_{\inc}$.\/}

\medskip\noindent
The following lemma is now proved as in \cite{st}.
\begin{lemma}\label{gooddifference}
Every good cell in $\boldsymbol{k}^n$ is of the form $X \setminus Y$
with $X,Y\in \operatorname{St}_n$.
\end{lemma}

\end{subsection}

\begin{subsection}{Good cell decomposition}  We obtain good cell
decomposition, namely, if $X_1 , \dots ,X_m \subseteq R^n$ are
definable, then there is a finite partition of $\boldsymbol{k}^n$
into good cells that partitions every $\st (X_i )$. A consequence of
this is that the $\boldsymbol{k}_{\ind}$-definable subsets of
$\boldsymbol{k}^n$ are finite unions of differences $\st
(X)\setminus \st (Y)$, where $X,Y \in \Def^n (R)$.

\medskip\noindent
\begin{lemma}\label{proj}
Let $C \subseteq \boldsymbol{k}^{n}$ be a good $i$-cell, let $X
\subseteq R^{n+1}$ be definable and suppose $k\in \{1,\dots,n\}$ is
such that $i_k=0$. Define $\pi \colon \boldsymbol{k}^{n+1}
\rightarrow \boldsymbol{k}^{n}$ by
$$\pi(x)= (x_1 , \dots , x_{k-1} , x_{k+1},
\dots ,x_{n+1} ).$$ Then $\pi \big(\st(X) \cap (C \times
\boldsymbol{k})\big)$ is a difference of sets in
$\operatorname{St}_n$.
\end{lemma}

\noindent
A {\em good decomposition of $I(\boldsymbol{k})^n$\/} is a special
kind of partition of $I(\boldsymbol{k})^n$ into finitely many good
cells. The definition is by recursion on $n$:
\begin{enumerate}
\item[(i)] a good decomposition of $I(\boldsymbol{k})$ is a collection $$\{(c_0 ,
c_1 ),(c_2 , c_3 ), \dots ,(c_k , c_{k+1}), \{c_0 \}, \{c_1 \} ,
\dots ,\{c_k \},\{c_{k+1} \} \}$$ of intervals and points in
$\boldsymbol{k}$ where $c_0 < c_1 < \dots < c_k < c_{k+1}$ are real
numbers with $c_0 = -1$ and $c_{k+1}=1$;
\item[(ii)] a good decomposition of
$I(\boldsymbol{k})^{n+1}$ is a finite partition $\mathcal{D}$ of
$I(\boldsymbol{k})^{n+1}$ into good cells such that
$\{\pi_{n}^{n+1}C:\; C \in \mathcal{D} \}$ is a good decomposition
of $I(\boldsymbol{k})^n$.
\end{enumerate}

\noindent

\begin{theorem}\label{thm} {\bf (Good Cell Decomposition)}
\begin{enumerate}
\item[{\em ($A_n$)}] Given any definable
$X_1 , \dots ,X_m \subseteq I^{n}$, there is a good decomposition of
$I(\boldsymbol{k} )^n$ partitioning each set $\st{X_i}$.

\item[{\em ($B_{n}$)}] If $f \colon X \rightarrow I$, with $X \subseteq
I^n$, is definable, then there is a good decomposition $\mathcal{D}$
of $I(\boldsymbol{k} )^n$ such that for every open $C \in
\mathcal{D}$, either the set $\st(\Gamma f) \cap (C \times
\boldsymbol{k})$ is empty, or $f$ induces a function $g\colon C
\rightarrow I(\boldsymbol{k} )$.
\end{enumerate}
\end{theorem}

\noindent Using the lemmas above the proof is very similar to that
of Theorem 4.3 in \cite{st}.

\medskip\noindent A {\em good decomposition of $\boldsymbol{k}^n$\/} is a
special kind of partition of $\boldsymbol{k}^n$ into finitely many
good cells. The definition is by recursion on $n$:
\begin{enumerate}
\item[(i)] a good decomposition of $\boldsymbol{k}^1=\boldsymbol{k}$
is a collection
$$\{(c_0 , c_1
),(c_2 , c_3 ), \dots ,(c_k , c_{k+1}), \{c_1 \} , \dots ,\{c_k \}
\}$$ of intervals and points in $\boldsymbol{k}$, where $ c_1 <
\dots < c_k \in \boldsymbol{k}$ and $c_0 = - \infty$, $c_{k+1}=
\infty$;
\item[(ii)] a good decomposition of $\boldsymbol{k}^{n+1}$ is a finite partition
$\mathcal{D}$ of $\boldsymbol{k}^{n+1}$ into good cells such that
$\{\pi^{n+1}_{n} C:\ C\in \mathcal{D}\}$ is a good decomposition of
$\boldsymbol{k}^n$.
\end{enumerate}

\begin{corollary}\label{reduction}
If $X_1 , \dots ,X_m \subseteq  R^n$ are definable, then there is a
good decomposition of $\boldsymbol{k}^n$ partitioning every
$\st{X_i}$.
\end{corollary}

\begin{theorem}\label{kinddefable}
The $\boldsymbol{k}_{\ind}$-definable subsets of $\boldsymbol{k}^n$
are exactly the sets of the form $\st (X) \setminus \st (Y)$ with
$X,Y \in \Def^n (R)$.
\end{theorem}

As in \cite{st} we obtain that the standard part of a
partial derivative of a definable function is almost everywhere
equal to the corresponding partial derivative of the standard part
of the function:


\begin{theorem}\label{derivatives}
Let $f \colon Y \rightarrow R$ with $Y \subseteq R^n$ be definable
with $V$-bounded graph. Then there is a good decomposition
$\mathcal{D}$ of $\boldsymbol{k}^n$ that partitions $\st{Y}$ such
that if $D \in \mathcal{D}$ is open and $D \subseteq \st{Y}$, then
$f$ is continuously differentiable on an open definable $X \subseteq
Y$ containing $D^h$, and $f, \frac{\partial f}{\partial x_1 }, \dots
, \frac{\partial f }{\partial x_n }$, as functions on $X$, induce
functions $g, g_1 , \dots , g_n \colon D \rightarrow \boldsymbol{k}$
such that $g$ is $C^1$ and $g_i = \frac{\partial g}{\partial x_i}$
for all $i$.
\end{theorem}

\end{subsection}

\end{section}

\begin{section}{The conditions $\mathcal{C}$, $\Sigma_{\inc}$, $\Sigma_{\dec}$ and $\Sigma$}
In this section we show that $\big(\Sigma_{\inc}\ \&\
\Sigma_{\dec}\big)$ implies $\mathcal{C}$, we prove that various
conditions imply $\Sigma$, and we give an example to the effect that
$\boldsymbol{k}_{\ind}$ is not always o-minimal.

\begin{subsection}{Closed and definably connected sets}
The conditions $\Sigma_{\dec}$ and $\mathcal{C}$ on pairs $(R,V)$
are stated on page \pageref{conditions}. Note that if $(R,V)$
satisfies $\mathcal{C}$, then $\boldsymbol{k}_{\ind}$ is o-minimal
by Lemma \ref{stan2}. For $(R,V)$ to satisfy $\mathcal{C}$ it
suffices that for each $n$ the closed
$\boldsymbol{k}_{\ind}$-definable subsets of $I(\boldsymbol{k})^n$
are exactly the sets $\st X$ with definable $X\subseteq I^n$. (This
follows by means of the homeomorphisms $\tau_n$.)


\begin{proposition}\label{closed}
Suppose $(R,V)\models \Sigma_{\inc}$ and $(R,V)\models
\Sigma_{\dec}$. Then $(R,V)$ satisfies $\mathcal{C}$. (In
particular, $\Sigma \Rightarrow \mathcal{C}$.)
\end{proposition}
\begin{proof}
The result will follow from Corollary \ref{reduction} once we show
that the closure of a good cell in $\boldsymbol{k}^{n}$ is of the
form $\st{X}$ for some definable $X \subseteq R^n$. Let $\epsilon$
range over all positive infinitesimals, and let $C \subseteq
\boldsymbol{k}^{n}$ be a good cell.

\medskip\noindent
{\em Claim.} There is $r_0 \in R^{>\ma}$ and a definable $X\subseteq
(0,r_0 )\times R^{n}$ such that
$$ 0 < r < r' < r_0\Longrightarrow X(r') \subseteq
X(r);   \qquad \st{\big(\bigcap_{\epsilon} X(\epsilon )\big)} = C
.$$

\smallskip\noindent
This claim follows by the same argument as the corresponding claim
in the proof of Proposition 5.1 in \cite{st}. Let $X\subseteq (0,r_0
)\times R^n$ be as in the Claim. Then, since $(R,V)\models
\Sigma_{\dec}$, we can take $\epsilon \in \ma^{>0}$ such that
$\st{X(\epsilon )} = \closure{C}$.
\end{proof}

For $Z\subseteq V^n$ we let $Z^h := \st^{-1}(\st (Z))$.
\begin{proposition}
Suppose $(R,V)$ satisfies $\mathcal{C}$, and let $X \subseteq V^n$
be definable and definably connected in $R$. Then $\st X$ is
definably connected.
\end{proposition}

\begin{proof}
Assume to the contrary that $\st{X}$ is not definably connected.
Then $\st{X}=\st{Y_1 } \dot\cup \st{Y_2 }$ for some definable,
nonempty $Y_1 , Y_2 \subseteq R^n$. We may assume that $Y_1 , Y_2$
are closed. Let
$$q:=
\inf\{d(y,\st Y_2 ):\; y \in \st{Y_1 } \}.$$ Since $\st{Y_1 } ,
\st{Y_2 }$ are closed and bounded, $q \in \k^{>0}$. Define
$$X_1 := \{x \in R^n :\; d(x,Y_1 ) \leq \frac{q}{4}  \}
\mbox{ and } X_2 := \{x \in R^n :\; d(x,Y_2 ) \leq \frac{q}{4} \}.$$
Then $X_1 , X_2$ are closed and disjoint, and $Y_{1}^{h} \subseteq
X_1$, $Y_{2}^{h} \subseteq X_2$. Since $X^h = Y_{1}^{h} \cup
Y_{2}^{h}$, we have $X = (X \cap X_1 ) \cup (X \cap X_2 )$, where $X
\cap X_1$, $X \cap X_2$ are nonempty, disjoint, and closed in $X$, a
contradiction with $X$ being definably connected.
\end{proof}

\end{subsection}

\begin{subsection}{Conditions implying $\Sigma$}

In the next lemma we use the following convention.  Let $C \subseteq
R^n$ be an $(i_1 , \dots ,i_n )$-cell of dimension $k$. Let
$$\lambda \colon \{ 1,\dots ,n \} \to \{1,\dots ,n  \}$$ be such that
$$1 \leq \lambda (1) < \dots < \lambda (k) \leq n$$ and
$i_{\lambda (1)} = \dots = i_{\lambda (k)}=1$. We define
$$C_0 := \{ a\in R^k :\ \mbox{ there is }x\in C \mbox{ such that }
x_{\lambda (1)}=a_{1}\, \& \dots \, \& \, x_{\lambda (k)}=a_{k}\}.$$
Then $C_0$ is the homeomorphic image of $C$ under a coordinate
projection $p \colon R^n \to R^k$. For a definable $C^1$-function
$f\colon C \to R$ we let $\hat{f}\colon C_0 \to R$ be defined by
$\hat{f}(p(x))=f(x)$ where $x \in C$. We denote by $\frac{\partial
f}{\partial x_j}(a)$, where $a\in C$ and $j\in \{ 1, \dots ,k\}$,
the $j$-th partial derivative of $\hat{f}$ at $p(a)$.

\begin{lemma}\label{cofinality}
Suppose $\operatorname{cofinality}(\m) > 2^{|\boldsymbol{k}|}$. Then
$(R,V)$ satisfies $\Sigma$.
\end{lemma}
\begin{proof}
Let $X\in \text{Def}^{1+n}(R)$. By cell decomposition we may assume
that $X$ is an $(i_1 , \dots , i_{n+1})$-cell satisfying for every
$k=1,\dots, n+1$ the following: If $p^{n+1}_{k} X = (f,g)$, then all
$\frac{\partial f}{\partial x_i }$, $\frac{\partial g}{\partial x_i
}$
have constant sign
on $p^{n+1}_{k-1}X$. If $p^{n+1}_{k}X=\Gamma f$, then all
$\frac{\partial f}{\partial x_i }$
have constant sign on $p^{n+1}_{k-1}X$.

Now there are $2^{|\boldsymbol{k}|}$ many distinct subsets of
$\boldsymbol{k}^n$. Let $f \colon \ma^{>0} \rightarrow \mathcal{P}(
\boldsymbol{k}^n )$, where $\mathcal{P}(\boldsymbol{k}^n)$ is the
power set of $\boldsymbol{k}^n$, be given by $\epsilon \mapsto \st
X(\epsilon )$. Assume to the contrary that for every $\epsilon_1 \in
\ma^{>0}$ we can find $\epsilon_2 \in \ma^{>\epsilon_1}$ such that
$\st X (\epsilon_1 )\not= \st X (\epsilon_2 )$. Then the above
assumption on $X$ yields a cofinal subset of $\ma$ such that $f$ is
injective on this subset, a contradiction.

\end{proof}
Note that, together with 5.3 and 6.4 in \cite{tconvex}, this lemma
implies that if $V$ is a $T$-convex subring of $R$, then
$(R,V)\models \Sigma$.

\begin{lemma}
Let $R$ be $\omega$-saturated. Then $(R,\mathcal{O})\models \Sigma$.
\end{lemma}
\begin{proof}
Let $X \subseteq R^{1+n}$ be defined over $a\in R^k$. Since $R$ is
$\omega$-saturated, we can take $\epsilon \in \ma$ such that
$\epsilon > \delta$ for every $\delta \in \text{dcl}(a)$ with
$\delta < \Q^{>0}$. Then for every $\epsilon'\in \ma^{>\epsilon}$,
$\text{tp}(\epsilon' |a)=\text{tp}(\epsilon |a)$, and, in
particular, $\st X(\epsilon' )=\st X (\epsilon )$. Otherwise we
could find $x\in \st X(\epsilon' ) \bigtriangleup \st X (\epsilon )$
and a box $B=(p_1 , q_1 )\times \dots \times (p_n , q_n )\subseteq
\R^n$ with $p_i , q_i \in \Q$ such that $x\in B$ and either
$\closure{B}\cap \st X(\epsilon )=\emptyset$ or $\closure{B}\cap \st
X(\epsilon' )=\emptyset$. Then $B' = (p_1 , q_1 )\times \dots \times
(p_n , q_n )\subseteq R^n$ is such that $B' \cap X(\epsilon
)=\emptyset$ and $B'\cap X(\epsilon') \not= \emptyset$, or vice
versa, a contradiction.
\end{proof}

We saw in Section 2 that if $(R,V)\models\Sigma_{\inc}$, then
$\boldsymbol{k}_{\ind}$ is o-minimal. However, the following example
shows that $\boldsymbol{k}_{\ind}$ is not always o-minimal.

\smallskip\noindent
{\bf Example. \/}\label{example} Let $\R_{\text{exp}}$ be the real
exponential field and let $R$ be a proper elementary extension. Take
$\lambda \in R$ such that $\lambda > \R$, and let $V$ be the
smallest convex subring of $R$ containing $\lambda$, i.e.
$$V:=\{y:\; |y|<\lambda^n \mbox{ for some }n \},$$ and let
$\boldsymbol{k}$ be the corresponding residue field. We define $\log
\colon R^{>0}\rightarrow R$ to be the inverse function of $\exp
\colon R \rightarrow R^{>0}$. Then $\log (V^{>0}) = V$ and it
induces an increasing and injective function
$\boldsymbol{k}^{>0}\rightarrow \boldsymbol{k}$, which, for
simplicity, we shall also denote by $\log$. Now the set $\{\st
(\lambda )^n : n\in \N \}$ is cofinal in $\boldsymbol{k}^{>0}$,
hence $\{\log{\st (\lambda )^n }:\; n\in \N \}$ is cofinal in
$\log{\boldsymbol{k}^{>0}}$. So the set $\log{\boldsymbol{k}^{>0}}$
is definable in $\boldsymbol{k}_{\ind}$, but, because
$\log{\st(\lambda )^n}=n\log{\st (\lambda )}$, it is not cofinal in
$\boldsymbol{k}^{>0}$, nor does it have a supremum. It follows that
$\boldsymbol{k}_{\ind}$ cannot be o-minimal, nor does $(R,V)$
satisfy $\Sigma_{\inc}$.

\end{subsection}

\end{section}

\begin{section}{Traces}

Recall from the Introduction that a set $X\subseteq R^n$ is a {\em
trace\/} if $X=Y\cap R^n$ for some $n$-ary relation $Y$ defined in
some elementary extension of $R$ using parameters from that
extension. Note that every $X\in \Def^n (R)$ is a trace, and that if
$X,Y \subseteq R^n$ are traces, then so are $X\cup Y$, $X\cap Y$ and
$X^{c}$. An example of a trace is $V\subseteq R$: take an element
$\lambda$ in an elementary extension of $R$ such that $V< \lambda <
R^{>V}$. Then $V=(-\lambda , \lambda)\cap R$ where the interval
$(-\lambda , \lambda)$ is taken in the extension.

We let $R^{\ast}$ be the expansion of $R$ by all traces $X\subseteq
R^n$, for all $n$. By the main result of \cite{bp} every subset of
$R^n$ definable in $R^{\ast}$ is a trace. It follows that every subset
of $R^n$ definable in $(R,V)$ is a trace.

\begin{lemma}\label{structure}
Let $\boldsymbol{k}^{\ast}$ be the expansion of the ordered field
$\boldsymbol{k}$ by the sets $\st ( X ) \subseteq \boldsymbol{k}^n$ for all
traces $X\subseteq R^n$ and all $n$. Then, for all $n$,
$$\Def^n (\boldsymbol{k}^{\ast})=\{ \st (X) :\; X\subseteq R^n
\mbox{ is a trace }\}.$$
\end{lemma}

\begin{proof}
We first show that for every $n$, the collection
$$\mathcal{C}_n := \{\st (X) :\; X \subseteq R^n
\mbox{ is a trace } \}$$ is a boolean algebra on
$\boldsymbol{k}^n$. It is clear that
$$\st( X_1 )
\cup \st ( X_2 )= \st ( X_1 \cup X_2 )$$ for all traces $X_1 , X_2
\subseteq R^n$. To see that $\mathcal{C}_n$ is closed under
complements, let $X\subseteq R^n$ be a trace, and note that
$$(\st X)^c = \st \{ y \in R^n :\; d(y , x) >
\ma \mbox{ for every }x\in X\}.$$ Since $\ma$ is a trace, the set
$\{ y\in R^n :\; d(y,x)
> \ma \mbox{ for all }x\in X\}$ is definable in $R^{\ast}$,
hence, by \cite{bp}, it
is itself a trace. We conclude that the sets $\st (X)$, where
$X\subseteq R^n$ is a trace, are the elements of a boolean algebra on
$\boldsymbol{k}^n$.

Now let $X\subseteq R^n$ be a trace, and let $0\leq m \leq n$. We
may assume that $X \subseteq V^n$ (since $V$ is a trace). Then
$\pi^{n}_{m}(\st (X)) = \st (p^{n}_{m} (X))$, and by \cite{bp},
$p^{n}_{m} ( X)$ is a trace.
\end{proof}

\noindent It follows from Lemma \ref{structure} that
$\boldsymbol{k}^{\ast}$ is weakly o-minimal.

\begin{lemma}\label{symdif} Let $S_1$ be a weakly o-minimal
structure and $S_2$ an o-minimal structure on the same underlying
ordered set $S$. Suppose for every $n$ and for every $X_1 \in
\Def^{n} (S_1 )$ there is $X_2 \in \Def^{n} (S_2 )$ such that $X_1
\bigtriangleup X_2$ has empty interior in $S^n$. Then
$\Def^{n}(S_1)\subseteq \Def^{n}(S_2 )$, for all $n$.
\end{lemma}
\begin{proof} We
proceed by induction on $n$. Let $n=1$. If $X\subseteq S$ is a
finite union of convex sets, and $Y \subseteq S$ is a finite union
of points and intervals, then either $X\bigtriangleup Y$ is finite,
or $X\bigtriangleup Y$ has nonempty interior. It follows that
$\Def^1 (S_1 )\subseteq \Def^1 (S_2 )$
 and, in particular, $S_1$ is o-minimal.

So assume $\Def^{k}(S_1 )\subseteq \Def^{k}(S_2 )$ holds for $k=1,
\dots ,n$. Since $S_1$ and $S_2$ are o-minimal, it suffices to show
that every $S_1$-cell in $S^{n+1}$ is definable in $S_2$. It is even
enough to prove this for $S_1$-cells $\Gamma g$; here $g\colon C \to
S$ is a continuous and $S_1$-definable function on an $S_1$-cell
$C\subseteq S^n$. Let  $\Gamma g$ be such an $S_1$-cell.

First, suppose $C$ is an open cell. By the inductive assumption $C
\in \Def^{n} (S_2 )$ and we can take $X \in \Def^{n+1} (S_2 )$ with
$X \subseteq C \times S$ such that $(-\infty , g) \bigtriangleup X$
does not contain a box. Let $p \colon S^{n+1} \to S^n$ be given by
$p (x_1,\dots, x_{n+1})=(x_1,...,x_n)$. For $X,Y\subseteq S^{n+1}$
we say that $X<Y$ if for all $a \in S^n$ and $(a,x)\in X$, $(a,y)\in
Y$ we have $x<y$. Now take an $S_2$-decomposition $\mathcal{D}$ of
$S^{n+1}$ which partitions $X$, and let $C_1 , \dots ,C_k$ be the
open cells in $p\mathcal{D}$ with $C_i \subseteq p X$. We claim that
$\Gamma (g|C_i ) \in \Def^{n+1}(S_2 )$ for every $i$.

So let $i\in\{1,\dots ,k \}$, and let $D_{1}, \dots ,D_{l}$ be the
open cells in $\mathcal{D}$ with $D_j \subseteq X$ and $p D_{j} =
C_i$ for all $j$. If $D_{j} = (f_j ,g_j )$ and $D_{j} \cap \Gamma (g
|C_i )\not= \emptyset$ for some $j\in \{1,\dots ,l \}$, then there
is $x \in C_i$ with $g(x)< g_j (x)$. Then, by continuity of $g$ and
$g_j$, we obtain a box $B \subseteq X \setminus (-\infty, g)$, a
contradiction. So $D_{j} \cap \Gamma g = \emptyset$, and, in
particular, $D_j < \Gamma (g|C_i )$ for every $j$.

Let $d \in \{ 1, \dots ,l \}$ be such that $D_j < D_{d} = (f_d , g_d
)$ for all $j\not=d$. If $g_d < g|C_i$ on a subset of $C_i$ with
nonempty interior, then, again by continuity of $g$ and $g_d$, we
find a box $B\subseteq (-\infty , g)$ with $\Gamma (g_d |p B) < B$.
Since $B$ intersects $X$ in only at most finitely many cells of the
form $\Gamma h$, where $h\colon C_i \rightarrow S$ is continuous, we
can find a box $B' \subseteq (-\infty , g)\setminus X$, a
contradiction. So $g_d = g|C_i$ outside a subset of $C_i$ with empty
interior, hence $g_d = g|C_i$ by continuity of $g$ and $g_d$.

We have shown that $\Gamma (g|C_i )$ is $S_2$-definable for all
$i=1, \dots ,k$. It is easy to check that then $$\Gamma g=\closure{
\bigcup_{i=1}^{k} \Gamma (g |C_i )} \cap (C_i \times S),$$ hence
$\Gamma g \in \Def^{n+1}(S_2 )$.

So let $\Gamma g \in \Def^{n+1}(S_2 )$ be an $(i_1 , \dots
,i_{n},0)$-cell with $i_k = 0$ where $1\leq k\leq n$, and let $$q
\colon S^{n+1} \rightarrow S^n \colon (x_{1}, \dots ,x_{n+1})
\mapsto (x_{1}, \dots ,x_{k-1}, x_{k+1}, \dots ,x_{n+1}).$$ By the
inductive assumption, $q (\Gamma g) \in \Def^{n}(S_2 )$. We define
$\Gamma g$ in $S_2$ as
$$\{ (x,y): \; x \in C \mbox{ and }(x_1 ,
\dots ,x_{k-1}, x_{k+1}, \dots ,x_{n},y)\in q(\Gamma g) \}.$$
\end{proof}

\medskip\noindent
The main result of this section is Theorem \ref{tracetheorem}, where
we assume that $R$ is $\omega$-saturated and $V=\mathcal{O}$. This
assumption is essential in that Theorem: Suppose
$\boldsymbol{k}_{\ind}$ is o-minimal but $\boldsymbol{k}$ is not
isomorphic to $\R$. Then $\boldsymbol{k}$ has a nonempty bounded
convex subset $X$ without a least upper bound in $\boldsymbol{k}$, so $X$ is not
definable in $\boldsymbol{k}_{\ind}$. However, $X^h \subseteq R$ is
a trace, and so $X=\st Y$ for some trace set $Y \subseteq R^n$.

\smallskip\noindent
{\bf In the rest of this section we assume that $R$ is
$\omega$-saturated and $V=\mathcal{O}$.\/} In particular, $\boldsymbol{k}=\R$.

\begin{lemma}\label{emptyint}
Let $Y\subseteq R^n$ be a trace. Then there is a definable $Z
\subseteq R^n$ such that $\st (Y) \bigtriangleup \st (Z)$ has empty
interior in $\R^n$.
\end{lemma}
\begin{proof} Take an elementary extension $R'$ of $R$ with a
definable set $Y'\subseteq R'^n$ such that $Y=Y'\cap R^n$. Then $Y'$
is defined in $R'$ by a formula $\phi(a,y)$ where $a\in R'^m$ and
$\phi(x,y)$ is a formula in the language of $R$, $x=(x_1,\dots,
x_m), y=(y_1,\dots, y_n)$. By $\omega$-saturation of $R$ we can take
$b\in R^m$ such that
$\text{tp}(b|\emptyset)=\text{tp}(a|\emptyset)$. Let $Z\subseteq
R^n$ be defined in $R$ by $\phi (b,y)$. Then $Y\cap\mathcal{O}^n
\subseteq \bigcup_{\epsilon } Z^{\epsilon }$, where $\epsilon$
ranges over all positive infinitesimals and $$Z^{\epsilon}:=\{ y\in
R^n :\; d(y,Z)\leq \epsilon \}.$$ Otherwise there would be $y\in
(Y\cap \mathcal{O}^n )$ such that $d(y,Z)>\ma$, so for some
$\mathcal{O}$-box $P\subseteq R^n$, we would have $P\cap Y\not=
\emptyset$ and $P\cap Z=\emptyset$, a contradiction with
$\text{tp}(b|\emptyset)=\text{tp}(a|\emptyset)$.

It follows that $\st(Y) \subseteq \st(Z)$. We claim that
$\interior{\st (Y) \bigtriangleup \st (Z)}=\emptyset$. Otherwise, we
can take a box $B\subseteq \R^n$ such that $B\subseteq \st (Z)
\setminus \st (Y)$, so the $V$-box lemma yields an $\mathcal{O}$-box
$P\subseteq Z$ such that $P\cap Y = \emptyset$, contradicting
$\text{tp}(b|\emptyset)=\text{tp}(a|\emptyset)$.
\end{proof}

\begin{theorem}\label{tracetheorem} For all $n$,
$$\Def^n (\R_{\ind })=\{ \st (X):\; X\subseteq R^n \mbox{ is a trace}
\}.$$
\end{theorem}
\begin{proof} By Lemma \ref{structure},
$$\{\st (X):\; X\subseteq R^n \mbox{ is a trace} \}=\Def^n
(\R^{\ast}),$$ for all $n$, and it is clear that $\Def^n (\R_{\ind
}) \subseteq \Def^{n}(\R^{\ast})$. So let $X\subseteq R^n$ be a
trace. By Lemma \ref{emptyint}, we can take $Y \in \Def^n (R)$ such
that $\interior{\st (X)\bigtriangleup \st Y} = \emptyset$, hence, by
Lemma \ref{symdif}, $\Def^n (\R^{\ast})\subseteq \Def^n
(\R_{\ind})$.
\end{proof}

\begin{corollary}\label{tracecor}
$\Def^n (\R_{\ind })=\{ \st (X):\; X\in \Def^n (R,\mathcal{O})\}$,
for all $n$.
\end{corollary}
\begin{proof}
It is clear that $\{ \st (X):\; X\in \Def^n (R,\mathcal{O})
\}\subseteq \Def^n (\R^{\ast})$, so by Theorem \ref{tracetheorem},
$\{ \st (X):\; X\in \Def^n (R,\mathcal{O}) \}\subseteq \Def^n
(\R_{\ind})$. To see that $$\Def^n (\R_{\ind}) \subseteq \{\st
(X):\; X\in \Def^n (R,\mathcal{O})  \},$$ recall that the
$\R_{\ind}$-definable subsets of $\R^n$ are finite unions of sets
$\st Y \setminus \st Z$, where $Y,Z\in \Def^n (R)$, and observe that
$$\st Y \setminus \st Z = \st \{ x\in Y :\; d(x,Z)>\ma \},$$
and that $\ma$ is definable in the structure $(R,\mathcal{O})$.
\end{proof}

\end{section}

\begin{section}{Open problems}
\begin{trivlist}
\item[1. ]We proved that $\Sigma_{\inc}$ implies o-minimality of
$\boldsymbol{k}_{\ind}$. Is the converse true?
\item[2. ]We showed that if $\text{cofinality}(\ma)> 2^{|\boldsymbol{k}|}$, then
$(R,V) \models\Sigma$. Conversely, if $(R,V)\models \Sigma$, is
there an elementary extension of $(R,V)$ satisfying this inequality?
\item[3. ]Does an analogue of Corollary \ref{tracecor} hold under more general
conditions, for example $(R,V)\models \Sigma$?
\item[4. ] Let $R$ be an $\omega$-saturated elementary extension of the
Lipshitz-Robinson structure. Are the definable sets of $\R_{\ind}$
just the semialgebraic sets?
\item[5. ]
The following question was posed by Lou van den Dries and Jonathan
Kirby:

\smallskip\noindent
$(\ast)$  Let $R$ be $\omega$-saturated and $V=\mathcal{O}$; is
$\R_{\ind}$ elementarily equivalent to a definable reduct of $R$?

\smallskip\noindent
To state this question precisely we assign to each $X\in \Def^n
(\R_{\ind})$ an $n$-ary relation symbol $P_X$, we let $L_{\ind}$ be
the language $L=\{<, 0,1, -,+, \cdot\}$ of ordered rings augmented
by these new relation symbols $P_X$, and we construe $\R_{\ind}$ as
a structure for the language $L_{\ind}$ in the obvious way, by
interpreting each $P_X$ as $X$. The formal statement of question
$(\ast)$ is as follows: does there exist an $L_{\ind}$-structure
$R'$ such that \begin{enumerate}
\item[(i)] $L$-reduct of $R'\ =\ L$-reduct of $R$,
\item[(ii)] each $n$-ary symbol $P_X$ is interpreted in $R'$ as a set
$X'\in  \Def^n (R)$,
\item[(iii)] $ \R_{\ind} \equiv R'$?
\end{enumerate}

\medskip\noindent
A positive solution might be hard to come by. To explain this, let
$L_{\exp}$ be the language of ordered rings augmented by a unary
function symbol $\exp$, and consider the $L_{\exp}$-theory
$T_{\exp}$ of the ordered exponential field $\R_{\exp}$. Peterzil
pointed out that by an argument as in Berarducci and Servi \cite{bs}
we have:

\begin{proposition} Suppose $(\ast)$ has a positive answer. Then
$T_{\exp}$ is decidable.
\end{proposition}
\begin{proof}
By \cite{bs} we have a recursive set $\Sigma_{\text{o}}$ of
$L_{\exp}$-sentences such that $T_{\text{exp}}\models \sigma$ for
all $\sigma\in \Sigma_{\text{o}}$ and all $L_{\exp}$-models of
$\Sigma_{\text{o}}$ are o-minimal. We can of course assume that
$\Sigma_{\text{o}}$ includes the usual axioms for real closed
fields, as well as an axiom expressing that $\exp$ is a
$C^1$-function with $\exp (0)=1$ and $\exp'=\exp $.

\smallskip\noindent
{\em Claim}. $\Sigma_{\text{o}}$ axiomatizes the (complete) theory
$T_{\exp}$.

\smallskip\noindent
To prove this claim, let $R$ be an $\omega$-saturated model of
$\Sigma_{\text{o}}$. Then the exponential function $\exp_R$ of $R$
induces the standard exponential function on $\R$. Since we assume
that $(\ast)$ has a positive answer for $R$, this gives a definable
function $e \colon R \to R$ such that $\R_{\exp}\equiv (R,e)$ (with
the last $R$ denoting its underlying ordered field). But this
function $e$ must be the exponential function $\exp_R$ by a
uniqueness result for solutions of differential equations in
o-minimal fields; see Otero, Peterzil and Pillay \cite{oka}. Thus
$\R_{\exp}\equiv R$.
\end{proof}

\end{trivlist}

\end{section}


\end{document}